\documentclass[a4paper, 12pt, reqno]{amsart}
\usepackage{indentfirst}
\usepackage{graphicx}
\usepackage{amsmath, amsthm}
\usepackage{amssymb}
\usepackage{txfonts}
\usepackage{comment}

\theoremstyle{definition}
\newtheorem{Def}{Definition}[section]

\theoremstyle{theorem}
\newtheorem{Thm}[Def]{Theorem}
\newtheorem{Coro}[Def]{Corollary}

\newtheorem{ThmZ}{Theorem}

\theoremstyle{definition}
\newtheorem{Remark}{Remark}[section]

  \makeatletter
  \@addtoreset{equation}{section}
  \makeatother

\def\C{\mathbb{C}}

\def\C{\mathbb{C}}

\title[Plurisubharmonicity of Bergman Kernels]{Plurisubharmonicity of Bergman Kernels on Generalized Annuli}
\author{Yanyan Wang}
\date{}

\thanks{
Project supported by NSFC No. 11031008 and NSFS No. 13ZR1444100.}

\address{Department of Mathematics
\\
Tongji University \\
Shanghai, 200092, China }
\email{wangyanyan0102@126.com}

\begin{document}

\maketitle {\bf Abstract:} {Let $A_\zeta=\Omega-\overline{\rho(\zeta)\cdot\Omega}$ be a family of generalized annuli over a domain $U$. We show that the logarithm $\log K_{\zeta}(z)$ of the Bergman kernel $K_{\zeta}(z)$ of $A_\zeta$ is plurisubharmonic provided $\rho\in PSH(U)$. It is remarkable that $A_\zeta$ is non-pseudoconvex when the dimension of $A_\zeta$ is larger than one. For standard annuli in ${\mathbb C}$, we obtain an interesting formula for  $\partial^2 \log K_{\zeta}/\partial \zeta\partial\bar{\zeta}$, as well as its boundary behavior.

}
\section{Introduction and Results}

\makeatother
\makeatletter
\renewenvironment{proof}[1][\proofname]{\par
  \pushQED{\qed}%
  \normalfont \topsep6\p@\@plus6\p@\relax
  \trivlist
  \item[\hskip\labelsep
        \textrm
    #1\@addpunct{.}]\ignorespaces
}{%
  \popQED\endtrivlist\@endpefalse
}
\makeatother
 \def\theequation{\thesection.\arabic{equation}}

\par In 2004, F. Maitani and H. Yamaguchi \cite{MY} brought a new viewpoint by studying the
variation of the Bergman metrics on the Riemann surfaces. Let us briefly recall their results.
\par Let $B$ be a disk in the complex $\zeta$-plane, $D$ be a domain in the product
space $B\times\mathbb{C}_{z}$, and let $\pi$ be the
first projection from $B\times\mathbb{C}_{z}$ to $B$ which is proper and smooth, and
$D_\zeta=\pi^{-1}(\zeta)$ a domain in $\C_z$. Let $K_\zeta$ denote the Bergman kernel of $D_\zeta$. Put $\partial
D=\bigcup_{\zeta\in B}(\zeta,\partial D_\zeta)$.

\begin{ThmZ}[cf. \cite{MY}] \label{maitani thm 1}
If $D$ is a pseudoconvex domain over $B\times\mathbb{C}_{z}$
with smooth boundary, then $\log
K_{\zeta}(z)$ is plurisubharmonic (psh) on $D$.
\end{ThmZ}
\begin{ThmZ}[cf. \cite{MY}]\label{maitani thm 2}
 If $D$ is a pseudoconvex domain over $B\times\mathbb{C}_{z}$
with smooth boundary, and for each $\zeta\in B$, $\partial
{D}$ has at least one strictly pseudoconvex point, then $\log K_{\zeta}(z)$ is a strictly psh
function on ${D}$.
\end{ThmZ}
In 2006, B. Berndtsson \cite{B.B2} made a striking generalization of Theorem \ref{maitani thm 1} to higher dimensional case, by using H\"ormander's $L^2-$estimates for $\bar{\partial}$:
\begin{ThmZ}[cf. \cite{B.B2}]\label{berndtsson thm}
 Let D be a pseudoconvex domain in $\C^k_\zeta\times \C^n_z$ and $\phi$ be a psh function on D. For each
$\zeta$ let $D_\zeta$ denote the $n$-dimensional slice
$D_\zeta:= \{z\in \C^n:\, (\zeta, z)\in D\}$
and by $\phi^\zeta$ the restriction of $\phi$ to $D_\zeta$. Let $K_{\zeta}(z)$ be the Bergman kernels of Bergman spaces $H^2(D_\zeta, e^{-\phi^\zeta})$. Then $\log K_{\zeta}(z)$ is psh
or identically equal to $-\infty$ on D.
\end{ThmZ}
\par These works rely heavily upon the pseudoconvexity of the total space ${D}$. In this paper, we obtain the plurisubharmonicity of $\log K_\zeta(z)$ for certain family of\/ {\it non-pseudoconvex} domains. Thus it would be interesting to find a more flexible condition than pseudoconvexity for the plurisubharmonicity of $\log K_\zeta(z)$.
     \par We consider the following family of generalized annuli
     $$
     A_\zeta=\Omega-\overline{\Omega_\zeta}
     $$
   where $\Omega\subset {\mathbb C}^n$ is a bounded complete circular domain and
   $$
   \Omega_\zeta=\rho(\zeta)\cdot \Omega:=\{\rho(\zeta) z:\, z\in \Omega\}
   $$
   with $0<\rho<1$ being an upper semicontinuous function on a domain $U$ in ${\mathbb C}^m$.
 Let $K_{\zeta}(z)$ denote the Bergman kernel of $A_\zeta$.

\begin{Thm}\label{annuli thm 3}
If $n\ge 2$ and $\rho\in PSH(U)$, then $\log K_{\zeta}(z)$ is a psh function on $U\times \Omega$. Furthermore, if $\rho$ is strictly psh on $U$, then $\log K_{\zeta}(z)$ is strictly psh on $U\times \Omega$.
\end{Thm}

The plurisubharmonicity of $\log K_\zeta(z)$ does not imply the pseudoconvexity of the total space even when the slices are planar domains. A simple example may be constructed as follows: let $D={\mathbb D}^2-\Gamma_f$ where $f$ is a\/ {\it non holomorphic} continuous self map of the unit disc ${\mathbb D}$ and $\Gamma_f$ is the graph of $f$. Since $\log K_\zeta(z)=\log K_{\mathbb D}(z)$, it is naturally psh, yet $D$ is not pseudoconvex, in view of a  theorem of Hartogs on holomorphicity of pseudoconcave continuous graphs. Nevertheless, it is still worthwhile to ask the following question.

\medskip

{\bf Question.} {\it Suppose $D$ is a bounded domain over ${U\times {\mathbb C}}$ where $U$ is a domain in ${\C}$. Let $K_\zeta$ denote the Bergman kernel of the slice $D_\zeta$. Suppose\/ $\log K_\zeta(z)$ is a psh function on $D$. Under which conditions is $D$ pseudoconvex?}

\medskip

  It is the case when $K_\zeta(z)\rightarrow \infty$ as $z\rightarrow \partial D_\zeta$ (note that $\log K_\zeta(z)$ is psh, in particular upper semicontinuous on $D$). We remind readers that Zwonek \cite{Zwonek} had given a complete characterization of Bergman exhaustiveness of bounded planar domains in terms of log capacities.

For standard annuli, i.e., $\Omega$ is the unit disc ${\mathbb D}$, $U$ is the punctured disc ${\mathbb D}^\ast$, and $\rho(\zeta)=|\zeta|$,  we have an interesting formula for $\partial^2\log K_\zeta/\partial\zeta\partial\bar{\zeta}$:
\begin{Thm}\label{annuli}
$$
\frac{\partial^2 \log K_{\zeta}(z)}{\partial \zeta\partial \overline{\zeta}}
=e^{2\omega_1}\frac{\left(2\mathcal{P}(u)-\mathcal{P}(\omega_1)+c\right)(\mathcal{P}(\omega_1)+c)}
{4\omega_1^{2}(\mathcal{P}(u)+c)^2}
$$
where $u=-2\log |z|$, $\omega_1=-\log|\zeta|$, $c(\omega_1)=\zeta(\omega_1)/\omega_1$, and $\mathcal{P}(\cdot)$ is the Weierstrass elliptic function function with periods  $2\omega_1,2\pi i$, and $\zeta(\cdot)$ is the Weierstrass zeta function.
\end{Thm}

As a consequence, we obtain

\begin{Coro}
 $\partial^2 \log K_{\zeta}(z)/\partial \zeta\partial
\overline{\zeta}\rightarrow 0$ as $D\ni (\zeta,z) \rightarrow \partial D$ in a nontrivial way, that is, at first $\zeta\rightarrow \zeta_0$, then $z\rightarrow \partial A_{\zeta_0}$.
\end{Coro}

\section{Proof of Theorem \ref{annuli thm 3}}
\begin{proof}
It is well-known that every holomorphic function $f$ on a bounded complete circular domain $\Omega$ admits a power series expansion as follows
$$f(z)=\sum_{j\geq0}p_j(z),$$
where $p_j(z)$ is a holomorphic polynomial of degree $j$, in the sense of locally uniform convergence. Thus the Bergman space $H^2(\Omega)$ of $\Omega$ admits a complete orthogonal basis
$$p_{j_1}, \cdots,p_{j_{m_j}}\in L_j,\,\,\,\,j=0,1,\cdots$$
where $L_j$ is the linear space spanned by homogeneous polynomials of degree $j$, and $m_j={\rm dim}_{\C}L_j$. Since
\begin{equation*}\label{}
  \int_{\Omega_\zeta}p_{j,r}\overline{p_{k,s}}=\rho(\zeta)^{2j+2k+2n}\int_{\Omega}p_{j,r}\overline{p_{k,s}}=0
\end{equation*}
for any pair $(j,r)\neq (k,s)$, it follows that
\begin{equation*}
 \int_{A_\zeta}p_{j,r}\overline{p_{k,s}}=\int_{\Omega}p_{j,r}\overline{p_{k,s}}-\int_{\Omega_\zeta}p_{j,r}\overline{p_{k,s}}=0.
\end{equation*}
By virtue of Hartogs' extension theorem, every holomorphic function on $A_\zeta$ can be extended to a holomorphic function on $\Omega$. Thus
\begin{equation}\label{}
   K_{\zeta}(z)=\sum_{j=0}^{\infty}\sum_{r=1}^{m_j}c_{j,r}\left|p_{j,r}(z)\right|^2
\end{equation}
where
\begin{eqnarray*}
c_{j,r}^{-1} &=& \int_{A_\zeta}\left|p_{j,r}(z)\right|^2= \int_{\Omega}\left|p_{j,r}(z)\right|^2- \int_{\Omega_\zeta}\left|p_{j,r}(z)\right|^2 \\
 &=& 1-\rho(\zeta)^{2j+2n}.
\end{eqnarray*}
That is,
\begin{equation}\label{}
  K_{\zeta}(z)=\sum_{j=0}^{\infty}\sum_{r=1}^{m_j}\frac{\left|p_{j,r}(z)\right|^2}{1-\rho(\zeta)^{2j+2n}}
\end{equation}
for any $z\in A_\zeta$. Put
$$
K_{\zeta}^{k}(z)=\sum_{j=0}^{k}\sum_{r=1}^{m_j}\frac{\left|p_{j,r}(z)\right|^2}{1-\rho(\zeta)^{2j+2n}}.
$$
Since $K^k_\zeta\in PSH(\Omega)$, we infer from the maximum principle that
$$
\max_{z\in M} K^k_\zeta(z)\le \max_{z\in\partial G} K^k_\zeta(z)\le  \max_{z\in\partial G} K_\zeta(z)
$$
where $M$ is a compact set whose interior contains $\overline{\Omega_\zeta}$ and $G$ a domain such that $M\subset G\subset\subset \Omega$. It follows immediately that the power series (2.2) converges uniformly on compact subsets of $\Omega$,  so that $K_\zeta$ can be extended to a smooth real function on $U\times \Omega$.
It is easy to verify that
\begin{equation*}\label{}
  u_{j}(\zeta,z)=\log \sum_{r=1}^{m_j}\left|p_{j,r}(z)\right|^2-\log\left(1-\rho(\zeta)^{2j+2n}\right )
\end{equation*}
is psh function on $\Omega$.
Since
\begin{equation}\label{}
  K_{\zeta}^{k}(z)=\sum_{j=0}^{k}e^{ u_{j}(\zeta,z)}
\end{equation}
and
$$
\chi(t_0,\cdots,t_k):=\log (e^{t_0}+\cdots+e^{t_k})
$$
 is a convex function which is non decreasing in each $t_j$,  we conclude that $\log K^k_\zeta(z)$ is psh on $U\times \Omega$ (see \cite{Demailly}, Theorem 4.16).
  Since $\left\{\log K^k_{\zeta}(z)\right\}^{\infty}_{k=0}$ is an increasing sequence of psh functions on $U\times \Omega$ whose limit is the\/ {\it continuous} function $\log K_{\zeta}(z) $, it follows that $\log K_{\zeta}(z) $ has to be psh on $U\times \Omega$.

  Now suppose $\rho$ is strictly psh on $U$. Without loss of generality, we may assume that the volume of $\Omega$ equals 1. Then
  $$
  u_0(\zeta,z)=u_0(\zeta)=-\log(1-\rho(\zeta)^{2n})
  $$
   is also strictly psh on $U$. Since $\chi$ is convex and non decreasing in each $t_j$, so
  $$
  \partial\bar{\partial}\log K^k_\zeta(z)\ge \frac{e^{u_0}}{K^k_\zeta(z)}\partial\bar{\partial} u_0(\zeta).
  $$
  Let $k\rightarrow \infty$, we get
   $$
  \partial\bar{\partial}\log K_\zeta(z)\ge \frac{e^{u_0}}{K_\zeta(z)}\partial\bar{\partial} u_0(\zeta),
  $$
 so that  for every $\xi=(\xi_1,\cdots,\xi_m,\xi_{m+1},\cdots,\xi_{m+n})$ with $(\xi_1,\cdots,\xi_m)\neq 0$, the Levi form $L(\log K_\zeta(z);\xi)>0$. While for every non zero vector $\xi=(0,\cdots,0,\xi_{m+1},\cdots,\xi_{m+n})$, we have
 $$
 L(\log K_\zeta(z);\xi)=\sum_{\alpha,\beta=1}^n \frac{\partial^2 \log K_\zeta(z)}{\partial z_j\partial\bar{z}_k}\xi_{m+\alpha}\overline{\xi_{m+\beta}}>0.
 $$
 Thus $\log K_\zeta(z)$ is strictly psh on $U\times \Omega$.
\end{proof}
 \medskip

 \begin{Remark} Since $\log K_\zeta(0)=u_0(\zeta)$, we conclude that $\log K_\zeta(z)$ would not be psh on $U\times \Omega$ if $u_0(\zeta)$ is not psh.
\end{Remark}
\section{Proof of Theorem \ref{annuli} and Corollary 1.3}

\begin{proof}[Proof of Theorem 1.2] It is known from \cite{NS} that
\begin{equation}\label{annuli 0}
\pi  K_{\zeta}(z)=\frac{\mathcal{P}(-2\log|z|)+\eta/(-\log|\zeta|)}{|z|^2},
\end{equation}
where
\begin{equation}\label{annuli eta}
 2\eta=\zeta(u-2\log|\zeta|)-\zeta(u),
\end{equation}
$u=-2\log |z|$,  $\mathcal{P}(\cdot)$ is the Weierstrass elliptic
function with periods $-2\log|\zeta|$, $2\pi i$, and $\zeta(\cdot)$
is the Weierstrass's zeta function.
If we let $\omega_1=-\log|\zeta|$, then (\ref{annuli 0}) changes to
\begin{equation}\label{annuli 2}
\pi  K_{\zeta}(z)=\frac{\mathcal{P}(u)+\eta/\omega_1}{|z|^2}.
\end{equation}
\par Since $\zeta'(\cdot)=-\mathcal{P}(\cdot)$, we have
$$\zeta'(\cdot+2\omega_1)=\zeta'(\cdot),$$
so that
$$\zeta(\cdot+2\omega_1)=\zeta(\cdot)+C.$$
Take $u=-\omega_1$, we get $C=2\zeta(\omega_1)$ and
\begin{equation}\label{annuli 3}
  \zeta(\cdot+2\omega_1)=\zeta(\cdot)+2\zeta(\omega_1).
\end{equation}
By (\ref{annuli eta}) and (\ref{annuli 3}), we  obtain
$\eta=\zeta(\omega_1)$.
Hence, (\ref{annuli 2}) changes to
\begin{equation}\label{annuli 1}
 K_{\zeta}(z)=\frac{\mathcal{P}(u)+c(\omega_1)}{\pi |z|^2},
\end{equation}
where $u= (0,2\omega_1)$,
$c(\omega_1)=\zeta(\omega_1)/\omega_1$.

 \par Now we going to calculate $\frac{\partial^2 }{\partial \zeta\partial \overline{\zeta}}\log K_{\zeta}(z)$. A straightforward calculation yields
 \begin{eqnarray*}
 \frac{\partial c(\omega_1)}{\partial\zeta}&=&\frac{\partial c(\omega_1)}{\partial\omega_1}\frac{\partial\omega_1}{\partial\zeta}
 =\frac{1}{2\zeta}\frac{\mathcal{P}(\omega_1)+c(\omega_1)}{\omega_1},\\
 \frac{\partial c(\omega_1)}{\partial\overline{\zeta}}&=&\frac{\partial c(\omega_1)}{\partial\omega_1}\frac{\partial\omega_1}{\partial\overline{\zeta}}
 =\frac{1}{2\overline{\zeta}}\frac{\mathcal{P}(\omega_1)+c(\omega_1)}{\omega_1},\\
 \frac{\partial^2 c(\omega_1)}{\partial\zeta\partial\overline{\zeta}}
 &=&\frac{\partial ^2 c(\omega_1)}{\partial\omega_1^2}\frac{\partial\omega_1}{\partial\zeta}\frac{\partial\omega_1}{\partial\overline{\zeta}}
 +\frac{\partial c(\omega_1)}{\partial\omega_1}\frac{\partial^2\omega_1}{\partial\zeta\partial\overline{\zeta}}\\
 &=&\frac{1}{4|\zeta|^2}\frac{\left(\mathcal{P}(\omega_1)+c(\omega_1)\right)-\omega_1\left(\mathcal{P}'(\omega_1)+c'(\omega_1)\right)}{\omega_1^2}.
 \end{eqnarray*}
 \par We claim that $\mathcal{P}'(\omega_1)=0$. To see this, simply note that $\mathcal{P}$ is an even function, hence $\mathcal{P}'(-\omega_1)=-\mathcal{P}'(\omega_1)$. Since
  $\mathcal{P}'(\omega_1)=\mathcal{P}'(-\omega_1)$ by periodicity, so we have $\mathcal{P}(\omega_1)=0$.
  It follows that
   $$\frac{\partial^2 c(\omega_1)}{\partial\zeta\partial\overline{\zeta}}=\frac{1}{4|\zeta|^2}\frac{2\left(\mathcal{P}(\omega_1)+c(\omega_1)\right)}{\omega_1^2},$$
So that
$$\frac{\partial^2 }{\partial \zeta\partial \overline{\zeta}}\log K_{\zeta}(z)
=e^{2\omega_1}\frac{\left(2\mathcal{P}(u)-\mathcal{P}(\omega_1)+c\right)(\mathcal{P}(\omega_1)+c)}
{4\omega_1^{2}(\mathcal{P}(u)+c)^2}.$$
\end{proof}

\begin{proof}[Proof of Corollary 1.3]
It is easy to see that  $\mathcal{P}(0)=\infty$ and  $\mathcal{P}(u)$ decreases in
$(0,\omega_1)$. We also know that
$\mathcal{P}(2\omega_1-u)=\mathcal{P}(u)$ and
$\omega_1^2\mathcal{P}(\omega_1)=\pi^2/6$. So
$\mathcal{P}(u)>0$ in $(0,2\omega_1)$.
 Note that
$$\mathcal{P}(u)=u^{-2}\left(1+O(u^2)\right)$$
as $u\rightarrow 0$.
 Thus,
$$2\mathcal{P}(u)-\mathcal{P}(\omega_1)+c=2u^{-2}\left(1+O(u^2)\right),$$
$$(\mathcal{P}(u)+c)^2=u^{-4}\left(1+O(u^2)\right).$$
If $|z|\rightarrow 1$, then $u
\rightarrow 0$.
Hence,
$$\lim_{|z|\rightarrow 1}\frac{\partial^2 \log
K_{\zeta}(z)}{\partial \zeta\partial \overline{\zeta}}=0.$$
 Using the periodicity of $\mathcal{P}(u)$,
we conclude that
$$\lim_{|z|\rightarrow |\zeta|}\frac{\partial^2 \log
K_{\zeta}(z)}{\partial \zeta\partial \overline{\zeta}}=0.$$

\end{proof}

\begin{Remark}
 The proof of Theorem \ref{annuli} implies that although the Levi form of $\log K_{D_\zeta}(z)$ with respect to $\zeta$ approaches to $0$ when $(\zeta,z)$ tends to the boundary of the domain,  $\log K_{D_\zeta}(z)$ is a strictly plurisubharmonic function on $D$. So, in Theorem \ref{maitani thm 2}, the condition that for each $\zeta\in B$,  $\partial
D$ has at least one strictly pseudoconvex point is only a sufficient condition for $\log K_{D_\zeta}(z,z)$ to be strictly psh  on $D$.
\end{Remark}
\begin{Remark}
The proof of Theorem 1.2 also yields the following equation
\begin{equation*}
  \frac{\partial^2 K_{\zeta}(z)}{\partial \zeta\partial \overline{\zeta}}=\frac{\partial K_{\zeta}(z) }{\partial \zeta}\frac{\partial K_{\zeta}(z)}{\partial \overline{\zeta}}.
\end{equation*}
\end{Remark}

{\bf Acknowledgements.}
The author is grateful for Professors Takeo Ohsawa and Bo-Yong Chen for their guidance.
 \vskip 4mm {

}


\begin{thebibliography}{10}

\bibitem {B.B2} B. Berndtsson, \textit{Subharmonicity properties of the Bergman kernel and some other functions
associated to pseudoconvex domains}, Ann. Inst. Fourier (Grenoble) \textbf{56} (2006), 1633-1662.
\bibitem {Demailly} J. P. Demally, Complex Analytic and Differential Geometry, http://www-fourier.ujf-grenoble.fr/~demailly/, 2012.
\bibitem {MY} F. Maitani and H. Yamaguchi, \textit{Variation of Bergman metrics on Riemann surfaces}, Math. Ann. \textbf{330} (2004), 477-489.
\bibitem {NS} N. Suita, \textit{Capacities and Kernels on Riemann Surfaces}, Arch. Ration. Mech. Anal. \textbf{46} (1972), 212-217.
\bibitem{Zwonek} W. Zwonek, {\it Wiener's type criterion for Bergman exhaustiveness}, Bull. Pol. Acad. Sci. Math. {\bf 50} (2002), 297--311.
\end{thebibliography}
\end{document}